\documentclass{amsart}
\usepackage[utf8]{inputenc}

\usepackage{amssymb,latexsym}

\usepackage{xcolor}
\usepackage{mathtools}
\usepackage[all,arc,curve,color,frame,pdf,line]{xy}
\usepackage{caption,comment}

\makeatletter
\renewcommand{\subsection}{\@startsection{subsection}{2}%
  {\z@}{1.5ex plus .5ex minus .2ex}{.5ex}%
  {\normalfont\bfseries}}
\makeatother

\newtheorem{theorem}{Theorem}[section]
\newtheorem{lemma}[theorem]{Lemma}
\newtheorem{proposition}[theorem]{Proposition}
\newtheorem{corollary}[theorem]{Corollary}
\newtheorem{problem}{Problem}
\newtheorem{definition}[theorem]{Definition}

\newtheorem{remark}[theorem]{Remark}
\newtheorem{example}[theorem]{Example}

\def\PSL{{\rm PSL}}
\def\PGammaL{{\rm P}\Gamma{\rm L}}

\def\soc{{\rm soc}}
\def\Z{{\mathbb Z}}

\def\la{{\langle}}
\def\ra{{\rangle}}

\def\fix{{\rm fix}}

\def\Aut{{\rm Aut}}

\def\gcd{{\rm gcd}}
\def\WStab{{\rm WStab}}
\def\Inn{{\rm Inn}}

\def\cal{\mathcal}

\newcommand{\ha}{\widehat a}
\newcommand{\hb}{\widehat b}
\newcommand{\hc}{\widehat c}
\newcommand{\hd}{\widehat d}
\newcommand{\hx}{\widehat x}
\newcommand{\hy}{\widehat y}
\newcommand{\hu}{\widehat u}

\usepackage[colorlinks]{hyperref}
\usepackage{cleveref}

\usepackage{tikz}
\usetikzlibrary{calc}

\title[Regular subgroups of wreath products]{Regular subgroups of wreath products and Cayley Cartesian products}

\author{Zhaochen Ding}
\address{Glagoljaška ulica 8, University of Primorska, 6000 Koper, Slovenia}
\email{zhaochen.ding@upr.si}

\author{Ted Dobson}
\email{ted.dobson@upr.si}

\author{Ademir Hujdurovi\'{c}}
\email{Ademir.Hujdurovic@upr.si}

\author{Yakun Wang}
\address{Department of Mathematics, Southern University of Science and Technology, Shenzhen 518055, P. R. China}
\email{wangyk3@sustech.edu.cn}

\date{\today}

\begin{document}

\begin{abstract}
The Cartesian product of Cayley graphs is again a Cayley graph.
For connected graphs, the converse question reduces to Cartesian
powers of Cartesian-prime graphs, and consequently to the study
of regular subgroups of wreath products in the product action.
We prove that, for a transitive permutation group $G$ of degree
$2p$ or $4p$, where $p$ is prime, $G\wr S_2$ contains a regular
subgroup if and only if $G$ does, and construct counterexamples
of degree $8p$ for every prime $p\equiv3\pmod4$.
By modifying this construction, we obtain an infinite family
of connected non-Cayley graphs whose Cartesian squares are
Cayley graphs, giving a negative answer to the converse question.
The smallest graph in our family has $120$ vertices and valency
$27$. In contrast, we prove that a connected vertex-transitive graph
of order $pq$, where $p$ and $q$ are primes, has a Cayley
Cartesian power if and only if the graph itself is Cayley.




\end{abstract}

\maketitle

\section{Introduction}

The relationship between graph products and their underlying components is one of the central themes in graph theory. A fundamental problem is to determine when a graph product inherits symmetry properties from its factors, and to what extent the symmetry of the product reflects the symmetry of those factors. 

Probably the most studied graph product is the Cartesian product. Recall that the \emph{Cartesian product}
$\Gamma_1\mathbin{\square}\cdots\mathbin{\square}\Gamma_k$ has vertex set
$V(\Gamma_1)\times\cdots\times V(\Gamma_k)$, with two vertices
adjacent if their entries are adjacent in one coordinate and equal
in all other coordinates. We write $\Gamma^{\mathbin{\square} n}$ for the
Cartesian product of $n$ copies of $\Gamma$.
The automorphism groups of Cartesian products are well
understood: They are determined by the automorphism groups of their
Cartesian-prime factors, together with permutations of isomorphic
factors \cite[Theorem 6.8]{HammackIK2011}. 

When studying symmetries of graphs, Cayley graphs are a central topic. By Sabidussi's theorem \cite{Sabidussi1958}, a graph is isomorphic
to a Cayley graph if and only if its automorphism group contains
a regular subgroup. Using the characterization of automorphism groups of Cartesian products,  it follows easily that a
Cartesian product of Cayley graphs is again a Cayley graph, as the direct product of regular subgroups of the automorphism groups of the factors acts regularly on the vertices of the product. 

The
converse question is considerably more difficult. Even with the
explicit description of the automorphism group, it is not clear
whether the existence of a regular subgroup forces the automorphism
groups of the factors to contain regular subgroups.
The structure of the automorphism group allows us to
reduce the question to a specific problem about permutation groups.
Two graphs are \emph{relatively prime} with respect to the Cartesian
product if they have no common nontrivial Cartesian factor.
If $\Gamma_1$ and $\Gamma_2$ are connected and relatively prime, then
$\Gamma_1\mathbin{\square}\Gamma_2$ is Cayley if and only if both factors
are Cayley; see \cite[Theorem 3.1]{HUJDUROVIC2018}.

By grouping isomorphic factors in the prime factorization of a
connected graph, the essential question therefore reduces to
Cartesian powers of connected Cartesian-prime graphs.
For such a graph $\Gamma$, we have
\[
\Aut(\Gamma^{\mathbin{\square} n})\cong\Aut(\Gamma)\wr S_n,
\]
where the wreath product acts in the product action on $V(\Gamma)^n$.
This leads to the following group-theoretic question: If $G$ is a
finite transitive permutation group and $G\wr S_n$ contains a
regular subgroup in the product action, does $G$ have to contain a regular
subgroup? 

The study of regular subgroups of finite primitive permutation
groups has a long history. Liebeck, Praeger, and Saxl
\cite[Corollary 3]{LiebeckPS2000} obtained a reduction theorem
that identified almost simple groups and groups in the product action
as the principal remaining cases in determining which primitive
groups contain regular subgroups. They subsequently classified
the regular subgroups of almost simple primitive groups
\cite{liebeck2010Regulara}, leaving the general product-action case
unresolved by those results. Another major contribution is
Li's classification of finite primitive permutation groups
containing an abelian regular subgroup \cite{Li2003}.
These results provide important structural information, but do
not settle the question of whether the existence of a regular
subgroup in $G\wr S_n$ forces a regular subgroup in $G$ for an arbitrary
transitive permutation group $G$.

An explicit example showing that the existence of a regular subgroup
in $G\wr S_n$ need not imply the existence of a regular subgroup in $G$
was given in \cite[Remark~3.7]{dobson2023When}: There is a transitive
group $G$ of degree $24$ and order $288$ with no regular subgroup
such that $G\wr S_2$, in the product action, contains a regular subgroup.

However, the graph question is more difficult than the group theoretic one, namely a negative answer to the group theoretic question does not automatically give
a graph counterexample. A graph admitting $G$ may have a larger
automorphism group containing a regular subgroup, even when $G$
itself contains none. The above-mentioned groups of degree 24 given in \cite{dobson2023When} did not produce an answer to the question about graphs.

In this paper, we give a negative
answer to the graph question by constructing an infinite family of connected non-Cayley
graphs whose Cartesian squares are Cayley graphs. To the best of
our knowledge, these are the first examples with this property.
Our construction gives the following result.

\begin{theorem}\label{thm:intro-counterexamples}
For every prime $p\equiv3\pmod4$, and odd integer $\ell\geq 3$ coprime to $p$ such that $-1$ is a square modulo
$\ell$, there exists a connected
non-Cayley graph $\Gamma$ of order $8p\ell$ and valency $8p+3$
such that $\Gamma\mathbin{\square} \Gamma$ is isomorphic to a Cayley graph.
\end{theorem}

The smallest graph in our constructed family has $120$ vertices
and valency $27$. To establish that these graphs are non-Cayley,
we determine their full automorphism groups and show that they
contain no regular subgroups, while their wreath products
with $S_2$ do.

In contrast, we prove that the converse holds for every Cartesian
power of a connected vertex-transitive graph whose order is a
product of two primes.

\begin{theorem}\label{thm:main}
Let $p$ and $q$ be primes, let $\Gamma$ be a connected
vertex-transitive graph of order $pq$, and let $n\ge1$ be an
integer. Then $\Gamma^{\mathbin{\square} n}$ is a Cayley graph if and only if $\Gamma$ is a Cayley graph.
\end{theorem}

The proof combines the classification of vertex-transitive graphs
of order $pq$ with restrictions on regular subgroups of wreath
products. In particular, it shows that no Cartesian power of the
Petersen graph is a Cayley graph, a fact for which we also give
a short direct proof.

We also prove several group-theoretic results. We show that if $G\wr S_n$ contains
a regular subgroup in the product action, then $G$ contains a
semiregular element of order $p$ for every prime $p$ dividing
its degree (\Cref{thm:semiregular-prime}). We also prove that
$G\wr S_n$ contains a regular subgroup if and only if $G$ does
whenever the degree of $G$ is relatively prime to $n!$
(\Cref{thm:coprime-n-factorial}).
For $n=2$, a further necessary condition is that $G$ contain
either a regular subgroup or a semiregular subgroup with exactly
two orbits. In particular, we prove that if $G$ has degree $2p$
or $4p$, where $p$ is prime, then $G\wr S_2$ contains a regular
subgroup in the product action if and only if $G$ contains a regular
subgroup (Proposition~\ref{prop:2p} and Theorem~\ref{thm:4p}).

The corresponding statement fails in degree $8p$:
For every prime $p\equiv3\pmod4$, Example~\ref{examp:8p} gives a
transitive group $G$ of degree $8p$ with no regular subgroup,
although $G\wr S_2$ contains one. These groups provide the
starting point for our graph construction. A modification
of the construction produces groups that can be realized
as the full automorphism groups of the required graphs.

The paper is organized as follows. In \Cref{sec:pre}, we recall
the automorphism structure of Cartesian products and establish
the reduction to wreath products. In \Cref{sec:general}, we
derive general restrictions on regular subgroups in product
action. We prove \Cref{thm:main} in \Cref{sec:pq-graph}.
In \Cref{sec:n=2}, we develop structural tools for wreath
products with $S_2$, which we apply in \Cref{sec:small} to
groups of degrees $2p$, $4p$, and $8p$. In the final section,
we construct the non-Cayley graphs with Cayley Cartesian
squares and determine their full automorphism groups.

\section{Preliminaries}\label{sec:pre}
We begin by recalling the structure of automorphism groups of Cartesian products.
Recall that two graphs are called \emph{relatively prime} if there 
exists no non-trivial graph that is a factor of both graphs with respect to the Cartesian product. A graph is said to be {\em prime} with respect to the Cartesian product if it cannot be factored as a Cartesian product of two non-trivial graphs. For a graph $X$, we denote the Cartesian product of $n$ copies of $X$ by $X^{\mathbin{\square} n}$.

\begin{theorem}\cite[Theorem 6.8]{HammackIK2011} 
\label{thm:automorphisms}
Let $\Gamma$ be a connected graph with prime factorization $\Gamma=\Gamma_1 \mathbin{\square} \Gamma_2 \mathbin{\square} \ldots \mathbin{\square} \Gamma_k$. Then for any automorphism $\varphi$ of $\Gamma$, there is a
permutation $\pi$ of $\{1, 2, \ldots , k\}$ and isomorphisms $\varphi_i : \Gamma_{\pi(i)} \to  \Gamma_i$ for which
$$
\varphi(x_1, x_2, \ldots , x_k)=(\varphi_1(x_{\pi(1)}),\ldots,\varphi_k(x_{\pi(k)})).
$$
\end{theorem}

Let $G_i\leq Sym(V_i)$ for $i\in \{1,\ldots,n\}$. The group $G_1\times G_2\times \ldots \times G_n$ acts {\em canonically} on $V_1\times V_2\times \ldots \times V_n$ in such a way that $g_i\in G_i$ is applied to the $i$-th coordinate.

The first observation explains why regular subgroups in a direct product force regular subgroups in each coordinate.

\begin{lemma}\label{lem:direct product}
Let $G_i\leq Sym(V_i)$ for $i\in \{1,\ldots,n\}$ be transitive groups. If there exists a regular subgroup of $G_1\times G_2\ldots \times G_n$ acting canonically on $V_1\times V_2\times \ldots \times V_n$, then every $G_i$ admits a regular subgroup.
\end{lemma}
\begin{proof}
Let $R$ be a regular subgroup of $G_1\times G_2\times\cdots\times G_n$. Since $R$ is regular, $|R|=|V_1|\cdots |V_n|$. Fix $j\in \{1,\ldots,n\}$, and choose vertices $v_i\in V_i$ for all $i\neq j$. Define
\[
S=\{(g_1,\ldots,g_n)\in R\mid g_i(v_i)=v_i\text{ for all }i\neq j\}
\]
and let
\[
S(j)=\{g_j\mid (g_1,\ldots,g_j,\ldots,g_n)\in S\}.
\]
Then $S(j)$ is a subgroup of $G_j$. Since $R$ is transitive on $V_1\times\cdots\times V_n$, the subgroup $S(j)$ is transitive on $V_j$. Moreover, if $S(j)$ were not semiregular, then $R$ would contain a non-identity element fixing a point of $V_1\times\cdots\times V_n$, contrary to the regularity of $R$. Thus $S(j)$ is regular on $V_j$. Since $j$ was arbitrary, the result follows.
\end{proof}

The following result follows directly by applying Lemma~\ref{lem:direct product} to the fact that the automorphism group of the Cartesian product of relatively prime graphs is the direct product of the automorphism groups of the factors (see \cite[Corollary 6.12]{HammackIK2011}; see also \cite[Theorem 3.1]{HUJDUROVIC2018}).

\begin{theorem}\label{thm:Cartesian_prime}
Let $X_1$ and $X_2$ be two connected relatively prime graphs with respect to the Cartesian product and let $X=X_1\mathbin{\square} X_2$. Then $X$ is isomorphic to a Cayley graph if and only if both $X_1$ and $X_2$ are Cayley graphs.
\end{theorem}

In light of \Cref{thm:Cartesian_prime}, the question of which Cartesian products are Cayley is reduced to the question: when is a Cartesian power of a graph isomorphic to a Cayley graph. The following result describes the automorphism group of a Cartesian power of a graph, which can be obtained directly from \Cref{thm:automorphisms}. 

\begin{lemma}
Let $\Gamma$ be a connected graph that is prime with respect to the Cartesian product. Then $\Aut(\Gamma^{\mathbin{\square} n})\cong \Aut(\Gamma) \wr S_n$.
\end{lemma}

\section{The general case for permutation groups}\label{sec:general}

Some of the preliminary results on regular subgroups and
Cartesian products were presented in the extended abstract of the Eurocomb 2023 conference \cite{dobson2023When}; we include their proofs here for completeness.


An element of a permutation group is said to be {\em semiregular} if all the cycles in its cyclic decomposition have equal lengths.

\begin{theorem}\label{thm:semiregular-prime}
Let $G$ be a transitive group of degree $d$, let $n\ge 2$ be an integer, and let $p$ be a prime divisor of $d$. If $G\wr S_n$ has a regular subgroup in the product action, then $G$ has a semiregular element of order $p$.
\end{theorem}

\begin{proof}

Let $R$ be a regular subgroup of $G\wr S_n$.
A Sylow $p$-subgroup of $R$ has order at least $p^n$. 
 On the other hand, by \cite[Lemma 3.1]{dobson2023When}, a Sylow $p$-subgroup of the coordinate-permuting group $S_n$ has order less than $p^n$. Hence, some element $\gamma\in R$ of order $p$ lies in the base group $G^n$. 

Write $\gamma=(\gamma_1,\ldots,\gamma_n)$. For each coordinate, $\gamma_i$ is either the identity or an element of order $p$ in $G$. Then, there is an element $\gamma_i$ of order $p$ without any fixed point. Otherwise, $\gamma$ is not semiregular, which contradicts that $\gamma\in R$.
\end{proof}

A direct corollary is as follows.
\begin{corollary}\label{Petersen tool}
    Let $\Gamma$ be a vertex-transitive digraph such that $\Aut(\Gamma)$ has no semiregular element of order $p$ for some prime $p$ dividing $|V(\Gamma)|$, and $\Gamma$ is prime with respect to the Cartesian product.  Then $\Gamma^{\mathbin{\square} n}$ is not isomorphic to a Cayley graph for every $n\ge 1$.
\end{corollary}
\begin{proof}
    We first observe that $\Gamma$ itself is not isomorphic to a Cayley graph. Indeed, if it were Cayley, then $\Aut(\Gamma)$ would contain a regular subgroup of order $|V(\Gamma)|$, and this subgroup would contain a semiregular element of order $p$, contrary to the hypothesis.

    Now suppose, for a contradiction, that $R$ is a regular subgroup of $\Aut(\Gamma^{\mathbin{\square} n})$. Since $\Gamma$ is prime with respect to the Cartesian product, \Cref{thm:automorphisms} gives $\Aut(\Gamma^{\mathbin{\square} n}) \cong \Aut(\Gamma)\wr S_n$ with the product action.
    This contradicts \Cref{thm:semiregular-prime}. This means that $\Aut(\Gamma^{\mathbin{\square} n})$ does not have any regular subgroup.
\end{proof}

\begin{corollary}\label{Petersen cor}
No Cartesian power of the Petersen graph is isomorphic to a Cayley graph.
\end{corollary}

\begin{proof}
By Corollary~\ref{Petersen tool}, it suffices to show that $\Aut(P)$ has no semiregular element of order $2$. The elements of order $2$ in $S_5$ are either transpositions or products of two transpositions. In the action of $S_5$ on the $2$-subsets of $\{1,2,3,4,5\}$, which is the automorphism group of the Petersen graph by \cite[Theorem 2.1.4]{Book}, every such involution fixes at least one $2$-subset. Hence no element of order $2$ in the automorphism group of the Petersen graph is semiregular.
\end{proof}

The next result gives a nice description of the existence of regular subgroups when the degree is relatively prime to $n!$.
\begin{theorem}\label{thm:coprime-n-factorial}
Let $G\leq Sym(V)$ be a transitive group whose degree is relatively prime to $n!$. If $G$ has no regular subgroup, then $G\wr S_n$ has no regular subgroup in the product action.
\end{theorem}
\begin{proof}
    Suppose that $R$ is a regular subgroup of $G\wr S_n$.
    Since $|R|=|V|^n$ and $|V|$ is relatively prime to $n!$, no nontrivial element of $R$ can induce a nontrivial permutation of the $n$ factors of $V^n$.
    Thus $R\le G^n$, and Lemma~\ref{lem:direct product} implies that $G$ has a regular subgroup, contrary to the hypothesis.
\end{proof}

The following two corollaries are straightforward.
\begin{corollary}
    Let $\Gamma$ be a vertex-transitive digraph that is not isomorphic to a Cayley digraph, whose automorphism group has order relatively prime to $n!$.  Then $\Gamma^{\mathbin{\square} n}$ is not isomorphic to a Cayley digraph.
\end{corollary}

\begin{corollary}
Let $\Gamma$ be a vertex-transitive digraph of odd order that is not isomorphic to a Cayley digraph and that is prime with respect to the Cartesian product. Then $\Gamma \mathbin{\square} \Gamma$ is not isomorphic to a Cayley digraph.
\end{corollary}

\section{Graphs of order $pq$}\label{sec:pq-graph}


In this section, we set $q<p$ to be two primes.
We start with an important observation about graphs of order $pq$ admiting a normal block system.

\begin{lemma}\label{lem:normal-block-system}
    Let $q<p$ be primes, and let $\Gamma$ be a vertex-transitive graph of order $pq$. If $\Aut(\Gamma)$ has a normal block system, then $\Gamma$ is isomorphic to a Cayley graph if and only if $\Aut(\Gamma)$ has a semiregular element of order $q$.
\end{lemma}
\begin{proof}
If $\Gamma$ is isomorphic to a Cayley graph, then $\Aut(\Gamma)$ contains a regular subgroup, and hence contains a semiregular element of order $q$. It remains to prove the converse. We argue by contradiction and assume that $\Aut(\Gamma)$ contains a semiregular element of order $q$, but that $\Gamma$ is not isomorphic to a Cayley graph. In particular, $\Gamma$ is not isomorphic to a circulant graph.

Let $\mathcal{B}$ be a normal block system with $r$ blocks, where $r=q$ or $p$. By \cite[Theorem 9.2.2]{Book},  it follows that $\Gamma$ is an $(r,pq/r)$-metacirculant. If $r=p$, then \cite[Exercise 1.5.6]{Book} implies that $\Gamma$ is a circulant graph of order $pq$, a contradiction.  Hence $r=q$, and the blocks of $\mathcal{B}$ have size $p$.

The automorphism groups of $(q,p)$-metacirculant graphs that are not isomorphic to circulant graphs whose automorphism groups have blocks of order $p$ (and the blocks are necessarily normal by \cite[Exercise 2.2.8]{Book}) are given in \cite[Theorem 3.2 (2)]{Dobson2006a}. Let $\alpha\in \Z_p^*$ and  define $\rho,\tau:\Z_q\times\Z_p\to\Z_q\times\Z_p$ by $\rho(i,j) = (i,j + 1)$ and $\tau(i,j) = (i + 1,\alpha j)$.  By \cite[Theorem 1.5.11]{Book}, we have that $\la\rho,\tau\ra\le\Aut(\Gamma)$ for some choice of $\alpha$.  In the statement of \cite[Theorem 3.2 (2)]{Dobson2006a}, it is assumed that $\alpha$ is {\it maximal} in $\Aut(\Gamma)$.  That is, whenever $\tau'(i,j) = (i + 1,\beta j)$ is contained in $\Aut(\Gamma)$, then $\la\rho,\tau'\ra\le\la\rho,\tau\ra$.  The existence of a maximal $\alpha$ is given by \cite[Lemma 2.7]{Dobson2006a}.   The second part (part $b$) of \cite[Theorem 3.2 (2)]{Dobson2006a} has two subcases, labeled i ($\vert\alpha\vert\not = q$ and $\Aut(\Gamma) = \la\rho,\tau\ra$) and ii (when $\vert\alpha\vert = q$ and so $\la\rho,\tau\ra$ is regular and $\Gamma$ is Cayley).


Suppose $\Aut(\Gamma) = \la\rho,\tau\ra$.  Write $|\alpha|=xy$, where $x$ is a power of $q$ and $y$ is relatively prime to $q$. Since $\Gamma$ is not isomorphic to a Cayley graph, \cite[Theorem 9]{AlspachP1982} gives $q^2\vert x$. By \cite[Lemma 1.5.1]{Book}, $|\Aut(\Gamma)|=|\la\rho,\tau\ra|=xyp$. Moreover, $\tau^y(i,j)=(i+y,\alpha^yj)$, and $\alpha^y$ has order $x$.  

Let $\gamma\in\Aut(\Gamma)$ be semiregular of order $q$. Then $\gamma$ lies in a Sylow $q$-subgroup of $\Aut(\Gamma)$. The subgroup $\la\tau^y\ra$ is such a Sylow subgroup, because it has order $x$ and $y$ is relatively prime to $q$. After conjugating if necessary, we may assume that $\gamma\in\la\tau^y\ra$. Write $x=q^b$, with $b\ge 2$. Since a cyclic group of order $q^b$ has a unique subgroup of order $q$, we have $1\not = \gamma\in\la(\tau^y)^{q^{b-1}}\ra$. However,
\[
(\tau^y)^{q^{b-1}}(i,j)=\tau^{xy/q}(i,j)=(i+xy/q,\alpha^{xy/q}j)=(i,\alpha^{xy/q}j).
\]
In particular, $\gamma(0,0)=(0,0)$, contradicting the assumption that $\gamma$ is semiregular. 


It only remains to consider the case \cite[Theorem 3.2 (2) (a)]{Dobson2006a}, where $\Aut(\Gamma)=\Z_2\ltimes F$ with
$F=\PSL(2,11)$ or $F=\PGammaL(n,k)$ with $n$ prime, $k=r^m$, $r$ prime, $\gcd(n,k-1)=1$, and $m$ a power of $n$. While the statement does not say so, the first paragraph of the proof of \cite[Theorem 3.2 (2)]{Dobson2006a} does; we have $q = 2$. Also, $F=\fix_{\Aut(\Gamma)}(\mathcal B)$ has blocks of size $p$.  There is an additional condition on the graphs in \cite[Theorem 3.2 (2) (a)]{Dobson2006a} that we will need but is not given in the statement, namely that $n\ge 3$.  In that proof, the automorphism group of $\Gamma$ was determined by first showing that $\Gamma$ is edge-transitive, and then using \cite[Table 1]{ChengO1987}.  Neither the statement of \cite[Theorem 3.2 (2)(a)]{Dobson2006a} nor \cite[Table 1]{ChengO1987} states that $n\ge3$, but this is a standing
assumption in \cite{ChengO1987} - see the start of the second paragraph on page 200, where the relevant graphs are defined.

Let $P=\la\rho\ra$ be a Sylow $p$-subgroup of $F$, which has order $p$, as $p = (k^n-1)/(k-1)$ is a primitive prime divisor of $k^n-1$ and $p \nmid (nm)$. By the Frattini argument, $\Aut(\Gamma)=F\,N_{\Aut(\Gamma)}(P)$, so $|N_{\Aut(\Gamma)}(P)|=2|N_F(P)|$. Now $|N_F(P)|$ is odd. If $F=\PSL(2,11)$, then $|N_F(P)|=55$. If $F=\PGammaL(n,k)$, then $P$ is a Singer cycle, so by \cite[II Satz 7.3]{Huppert1967}, $|N_F(P)|=pnm$, which is odd because $p>q=2$, $n\ge3$ is an odd prime, and $m$ is a power of $n$.  Hence a Sylow $2$-subgroup of $N_{\Aut(\Gamma)}(P)$ has order $2$ and is generated by an involution $\sigma\notin F$. Then $\sigma$ interchanges the two blocks, and $\la P,\sigma\ra$ is a transitive subgroup of order $2p$, hence regular. Thus $\Gamma$ is a Cayley graph, a contradiction.  
\end{proof}

\begin{remark}
     For $q>2$, if $\Aut(\Gamma)$ does not have a normal block system, then the above theorem is not necessarily true.
    For example, the line graph of the Petersen graph has $15$ vertices, and is not isomorphic to a Cayley graph, but its automorphism group contains semiregular automorphisms of order $3$.
\end{remark}

We now divide the analysis of vertex-transitive graphs of order $pq$ into two cases, according to whether $q=2$.

\begin{theorem}\label{thm:metacirculant}
Let $p>2$ be a prime and let $\Gamma$ be a vertex-transitive graph of order $2p$. Then $\Gamma$ is isomorphic to a Cayley graph if and only if $\Aut(\Gamma)$ has a semiregular element of order $2$.
\end{theorem}

\begin{proof}
If $\Gamma$ is a Cayley graph, then $\Aut(\Gamma)$ contains a regular subgroup $R$, and so a Sylow $2$-subgroup of $R$ has order $2$ and is semiregular.

Conversely, first suppose that $\Aut(\Gamma)$ is primitive. By \cite{LiebeckS1985a}, the primitive groups of degree $2p$ are either $2$-transitive, or else $p=5$ and the group is the action of $A_5$ or $S_5$ on the $2$-subsets of $\{1,2,3,4,5\}$. In the first case, $\Gamma$ is complete or has no edges, so $\Aut(\Gamma)=S_{2p}$.  Thus $\Gamma$ is isomorphic to a Cayley graph. In the second case, the only graphs with this automorphism group are the Petersen graph and its complement, and we saw in Corollary \ref{Petersen cor} that their common automorphism group contains no semiregular element of order $2$.


Suppose $\Aut(\Gamma)$ is imprimitive.  In view of Lemma \ref{lem:normal-block-system}, we may assume $\Aut(\Gamma)$ does not have a normal block system.  As a block system with blocks of size $p$ is always normal \cite[Exercise 2.2.8]{Book}, the only remaining possibility is that $\Aut(\Gamma)$ has only block systems with blocks of size $2$, and none of them is a normal block system. By \cite[Theorem 3.2 (1)]{Dobson2006a}, all such graphs have odd order, so this case cannot occur for graphs of order $2p$. 
\end{proof}

\begin{theorem}\label{thm:not-metacirculant}
    Let $q<p$ be odd primes, and let $\Gamma$ be a vertex-transitive graph of order $pq$.
    If $\Gamma$ is not metacirculant, then $\Gamma^{\mathbin{\square} n}$ is not isomorphic to a Cayley graph for every $n\ge 1$.
\end{theorem}
\begin{proof}
    We argue by contradiction and assume that $\Gamma^{\mathbin{\square} n}$ is isomorphic to a Cayley graph for some $n\ge 1$.
    If $n=1$, then \cite[Proposition 1.2]{MarusicS1992} would imply that $\Gamma$ is metacirculant, contrary to the hypothesis. Hence $n\ge 2$. Also observe that $\Gamma$ must be connected; otherwise, $\Gamma$ is a disjoint union of graphs of prime order, implying that $\Gamma$ is circulant, contrary to the assumption that it is not metacirculant.

    We first reduce to the case where $\Gamma$ is prime with respect to the Cartesian product. If not, then $\Gamma=\Sigma_1\mathbin{\square}\Sigma_2$, where $\Sigma_1$ has order $p$ and $\Sigma_2$ has order $q$.
    Since $\Gamma$ is vertex-transitive, both $\Sigma_1$ and $\Sigma_2$ are vertex-transitive by \cite[Theorem 6.17]{HammackIK2011}.
    Thus both factors are circulant, and so $\Gamma$ is isomorphic to a Cayley graph. This would again make $\Gamma$ metacirculant, a contradiction. Therefore $\Gamma$ is prime with respect to the Cartesian product.

    It follows that $\Aut(\Gamma^{\mathbin{\square} n})=\Aut(\Gamma)\wr S_n$. By Sabidussi's theorem, there is a regular subgroup $R\le\Aut(\Gamma)\wr S_n$. In particular, $|R|=p^nq^n$. Burnside's theorem implies that $R$ is solvable.
    Let $\pi:\Aut(\Gamma)\wr S_n\rightarrow S_n$ be the natural projection and let $p_1:\Aut(\Gamma)^n\rightarrow\Aut(\Gamma)$ be the natural projection onto the first coordinate.
    If $\pi(R)$ is transitive, then by \cite[Proposition 2.5]{LiebeckPS2000}, $p_1(R\cap \Aut(\Gamma)^{n})$ is a transitive subgroup of $\Aut(\Gamma)$.
    Suppose that $\pi(R)$ is intransitive and has an orbit of length $k<n$.
    Then $R\le\Aut(\Gamma)\wr S_{k}\times\Aut(\Gamma)\wr S_{n-k}$.
    Let $q_1:\Aut(\Gamma)\wr S_{k}\times\Aut(\Gamma)\wr S_{n-k}\rightarrow\Aut(\Gamma)\wr S_{k}$ be the natural projection.
    By \cite[Proposition 2.4 {\rm (ii)}]{LiebeckPS2000}, $q_1(R)$ is a transitive subgroup of $\Aut(\Gamma)\wr S_{k}$ on $(V\Gamma)^k$.
    Let $q_2:\Aut(\Gamma)^k\rightarrow\Aut(\Gamma)$ be the natural projection.
    Applying \cite[Proposition 2.5]{LiebeckPS2000} gives that $q_2(q_1(R)\cap (\Aut(\Gamma))^k)$ is a transitive subgroup of $\Aut(\Gamma)$.
    Therefore, $\Aut(\Gamma)$ always has a transitive $\{p,q\}$-subgroup.
    Denote this subgroup by $L$.

   Since \(L\) is solvable and has degree \(pq\), it cannot be primitive; that is, a finite solvable primitive permutation group has prime-power degree. Hence \(L\) admits a nontrivial block system, whose blocks have size \(p\) or \(q\).
If the blocks are of size $p$, then by \cite[Exercise 2.2.8]{Book}, this is a normal block system, which implies that $\Gamma$ is metacirculant by \cite[Theorem 9.2.2]{Book}, contrary to the hypothesis.
Therefore the blocks have size \(q\). 
Since \(\Gamma\) is not metacirculant, by \cite[Theorem 3.4]{Marusic1988},  \(L\) is nonsolvable, a contradiction.

    The contradiction shows that no such regular subgroup $R$ exists. Hence $\Gamma^{\mathbin{\square} n}$ is not isomorphic to a Cayley graph for every $n\ge 1$.
\end{proof}

\begin{proof}[Proof of Theorem~\ref{thm:main}]
    Let $q\leq p$ be primes with $|V(\Gamma)|=pq$, and assume that $\Gamma$ is not isomorphic to a Cayley graph. By \cite{Marusic1985}, every vertex-transitive graph of order $p^2$ is Cayley, hence we may assume that $q<p$.
    If $q=2$, the conclusion follows from Theorem~\ref{thm:metacirculant} and Corollary~\ref{Petersen tool}.
    We may therefore assume that $q>2$.

    By Theorem~\ref{thm:not-metacirculant}, it remains only to treat the case where $\Gamma$ is metacirculant. 
    If $\Aut(\Gamma)$ has a normal block system,
    Lemma~\ref{lem:normal-block-system} then shows that $\Aut(\Gamma)$ has no semiregular element of order $q$ and
    Corollary~\ref{Petersen tool} now applies with this prime $q$, and gives that $\Gamma^{\mathbin{\square} n}$ is not isomorphic to a Cayley graph for every $n\ge 1$.
    Suppose now that $\Aut(\Gamma)$ does not have a normal block system.
    By \cite[Exercise 2.2.8]{Book}, either $\Aut(\Gamma)$ has a block of size $q$ or $\Aut(\Gamma)$ is primitive.
    For the former case, \cite[Theorem 12]{DobsonHKM2025} implies that $\Gamma$ is isomorphic to a Cayley graph or Maru\v si\v c-Scapellato graph that is not isomorphic to a metacirculant graph, a contradiction. We conclude that $\Aut(\Gamma)$ is primitive.
    By \cite[pp. 192]{MarusicS1994}, $\soc(\Aut(\Gamma))=\PSL(2,19)$ and its action on $V\Gamma$ is permutation equivalent to the coset action of $\PSL(2,19)$ on cosets of $A_5$. 
    As $\soc(\Aut(\Gamma))=\PSL(2,19)$, $\Aut(\Gamma)=\PSL(2,19)$ or $\operatorname{PGL}(2,19)$.
    For each case, $\Aut(\Gamma)$ has cyclic Sylow $3$-subgroups of order $9$. Hence, all subgroups of order $3$ in $\Aut(\Gamma)$ are conjugate.
    Since $A_5$ has a Sylow $3$-subgroup of order $3$, it follows that every subgroup of $\Aut(\Gamma)$ of order $3$ is conjugate to a subgroup of a point stabilizer, which means that $\Aut(\Gamma)$ does not have a semiregular element of order $3$.
    By Corollary~\ref{Petersen tool}, $\Gamma^{\mathbin{\square} n}$ is not isomorphic to a Cayley graph for every $n$.
    This completes the proof.
\end{proof}

\section{Wreath product by the Symmetric group of degree 2}\label{sec:n=2}
In this section, we study when $G\wr S_2$ can contain a regular subgroup in the product action for a permutation group $G$. 
We first present a general necessary condition and then provide some preliminary results, which will be used in \Cref{sec:small}.
\begin{theorem}
Let $G\leq Sym(V)$ be transitive. If $G\wr S_2$ contains a regular subgroup (in the product action), then $G$ admits a regular subgroup or $G$ admits a semiregular subgroup with two orbits.
\end{theorem}
\begin{proof}
Let $R\leq G\wr S_2$ be a regular subgroup and let $d=|V|$. If $R\leq G\times G$, then Lemma~\ref{lem:direct product} implies that $G$ admits a regular subgroup. We may therefore assume that $R$ is not contained in $G\times G$. Let $B=R\cap (G\times G)$. Then $B$ has index two in $R$, and since $R$ is regular, $B$ is semiregular with two orbits on $V\times V$, hence $|B|=\frac{d^2}{2}$.

 Consider the action of $B$ on $V$ through the second coordinate, that is, $(b_1,b_2)\in B$ sends $x\in V$ to $b_2(x)$. Let $v\in V$, and let $K=\{(b_1,b_2)\in B\mid b_2(v)=v\}$. Observe that $K$ is the stabilizer of $v$ in the action of $B$ on $V$. By the orbit stabilizer theorem, it follows that $|B|=|K|\cdot |v^B|\leq |K|d$. Since $|B|=\frac{d^2}{2}$, it follows that $|K|\geq \frac{d}{2}$. 

 Consider the projection $\pi_1:K\to G$ defined by $\pi_1(b_1,b_2)=b_1$. Let $H=\pi_1(K)\leq G$. We claim that $H$ is semiregular with at most two orbits. Since $K\leq R$ and $R$ is regular, it follows that $H$ is semiregular. Suppose that $(b_1,b_2)\in K$ belongs to the kernel of $\pi_1$. Then $b_1=1$ and $b_2(v)=v$, which implies that $(b_1,b_2)$ fixes $(v,v)$. Since $(b_1,b_2)\in B\leq R$, and $R$ is regular, it follows that $(b_1,b_2)=(1,1)$, hence $\pi_1$ is injective. It follows that $|H|=|K|\geq \frac{d}{2}$. Since $H$ is semiregular, it now follows that either $H$ has two orbits of size $d/2$, or it is regular.  
 \end{proof}


\begin{corollary}
Let $\Gamma$ be a connected non-Cayley digraph that is prime with respect to the Cartesian product such that $\Aut(\Gamma)$ admits no semiregular subgroup with two orbits. Then $\Gamma\mathbin{\square} \Gamma$ is not isomorphic to a Cayley digraph.
\end{corollary}

We collect several elementary tools used in \Cref{sec:small}.

\begin{lemma}\label{lem:regular-quotient-something}
    Let $G$ be a group, let $H$ be a subgroup and let $N$ be a normal subgroup such that $H\cap N=1$.
    \begin{itemize}
        \item [\rm (i)] If $HN/N$ has a complement in $G/N$, then $H$ has a complement in $G$.
        \item [\rm (ii)] If $H$ has a complement containing $N$ in $G$, then $HN/N$ has a complement in $G/N$.
        \item [\rm (iii)] If $\gcd(|N|,|G/N|)=1$, then $H$ has a complement in $G$
if and only if $HN/N$ has a complement in $G/N$.
    \end{itemize}
    
\end{lemma}
\begin{proof}
    {\rm (i)} Let $R/N$ be a complement of $HN/N$ in $G/N$.
    We claim that $R$ is a complement of $H$.
    Since $(R/N)\cap (HN)/N=1$, we have $R\cap H\le R\cap HN\le N$.
    Therefore, $R\cap H\le N\cap H=1$.
    Also note that $G=RNH=RH$. This concludes {\rm (i)}.

    {\rm (ii)} Let $R$ be a complement of $H$ that contains $N$.
    Then $G/N=(R/N)(HN/N)$ and $(R/N)\cap (HN)/N=(R\cap HN)/N=N/N=1$.
    Therefore, $R/N$ is a complement of $HN/N$ in $G/N$.

    {\rm (iii)} Let $R$ be a complement of $H$ in $G$.
    We claim that $N\le R$. 
    Since $N\cap H=1$ and $\gcd(|N|,|G/N|)=1$, $\gcd(|N|,|H|)=1$.
    Then by $|G|=|R||H|$, $|N|\mid |R|$.
    Since $RN\le G$, we have $(|R||N|/|R\cap N|)\mid |G|$.
    Then $|N|/|R\cap N|\mid |H|$.
    As $\gcd(|N|,|H|)=1$, $|N|=|R\cap N|$, which implies that $N\le R$.
    Then {\rm (ii)} implies that $R/N$ is a complement of $HN/N$ in $G/N$. 
\end{proof}

Let $\mathcal{K}_2(d)$ be the class of transitive groups $G$ of degree $d$ such that $G$ does not have a regular subgroup, but $G\wr S_2$ with the product action has a regular subgroup. 
We next study minimal groups in $\mathcal{K}_2(d)$ with respect to inclusion. This will be useful for proving that $\mathcal{K}_2(d)$ is empty for certain values of $d$.

\begin{lemma}\label{lem:reduction}
Let $G\in\mathcal{K}_2(d)$, and suppose that no proper transitive subgroup of $G$
belongs to $\mathcal{K}_2(d)$. Let $H=G_\alpha$ and $h=|H|$. If $R$ is regular in
$G\wr S_2$, then, after conjugating $R$ inside $G\wr S_2$,
$B=R\cap G^2$ is a subdirect product in $G^2$. There is a semiregular normal subgroup
$N\trianglelefteq G$ such that
\[
 |B|=\frac{d^2}{2},\qquad
 |N|=\frac{d}{2h},\qquad
 |G/N|=2h^2.
\]
Moreover, $B/N^2\le (G/N)^2$ is the graph of an automorphism $f$ of $G/N$, and
$HN/N\cong H$. Denoting $A=HN/N$, we have $A\cap f(A)=1$ and  $A$ has no complement in $G/N$.
\end{lemma}

\begin{proof}
By Lemma~\ref{lem:direct product}, $R\not\le G^2$.
Thus $|R:B|=2$ and $|B|=\frac{d^2}{2}$. Conjugating in the base group, we may choose
$\tau=(g,1)\sigma\in R\setminus B$, where $\sigma$ swaps the coordinates.
Since $\tau^2=(g,g)\in B$ and $\tau$ normalizes $B$, the two projections of $B$
coincide. Let $S\leq G$ be the projection. Then $R\leq S\wr S_2$ and $S\wr S_2$ is transitive in the product action. The orbit of
$(\alpha,\alpha)$ under $S\wr S_2$ is $\alpha^S\times \alpha^S$, so $S$ is transitive.
Minimality gives $S=G$.

Let $N_1=\{x\in G\mid (x,1)\in B\}$ and $N_2=\{y\in G\mid (1,y)\in B\}$.
By Goursat's lemma, both groups are normal in $G$. Identify $G\times 1$ and $1\times G$ with $G$. Then conjugation by $\tau$ gives
$N_1=N_2$. Write $N_1=N_2=N\le G$. Since $N\times1\leq B$, the group $N$ is a semiregular subgroup of $G$.
The projection of $N$ onto the first coordinate gives $B/(1\times N_2)\cong G$, so $|B|=|G||N|$.
Combining this with $|G|/h=d$ and $|B|=\frac{d^2}{2}$, we obtain $|N|=\frac{d}{2h}$ and $|G/N|=2h^2$.
Modulo $N^2$, the two projections are isomorphisms, giving the automorphism $f$.

Suppose that $A\cap f(A)\ne1$, and let $1\neq x\in A\cap f(A)$. Then $x=h_2N$ for some $1\neq h_2\in H$, and also $x=f(h_1N)$ for some $1\neq h_1 \in H$. The element $(h_1N,h_2N)=(h_1N,f(h_1N))$ then belongs to $B/N^2$, and hence $(h_1,h_2)\in B$. This is a contradiction to $B\cap (H\times H)=1$.
If $A$ has a complement $T/N$ in $G/N$, then $T$ is a complement of $H=G_\alpha$ in $G$ by Lemma~\ref{lem:regular-quotient-something}, hence $T$ is regular.
\end{proof}

\begin{corollary}\label{cor:minimal-counterexample}
    Let $G\in \mathcal{K}_2(d)$ be a transitive group such that no proper transitive subgroup of $G$ belongs to $\mathcal{K}_2(d)$.
    If $G$ has an element of order $p^n$ for a prime $p$ and $n\ge 1$, then $G$ has a semiregular element of order $p^n$.
    In particular, $p^n\mid d$.
\end{corollary}
\begin{proof}
    Let $x\in G$ have order $p^n$.
    By Lemma~\ref{lem:reduction}, $B$ has an element of order divisible by $p^n$.
    Since $B$ is semiregular, $B$ has a semiregular element $b$ of order $p^n$.
    Write $b=(u,v)$. Its power $b^{p^{n-1}}$ is fixed-point-free, so one of $u^{p^{n-1}}$ and $v^{p^{n-1}}$ is fixed-point-free. 
    Without loss of generality, assume $u^{p^{n-1}}$ is fixed-point-free. Then every cycle of $u$ has length $p^n$. This completes the proof.
\end{proof}

\begin{corollary}\label{cor:prime-stabilizer}
Let $G$ be a transitive group with point stabilizer of prime order. Then $G\wr S_2$ has a regular subgroup (in the product action) if and only if $G$ has a regular subgroup.
\end{corollary}
\begin{proof} Let $p=|G_\alpha|$ and let $d$ denote the degree of $G$. If $G$ has a regular subgroup $H$, then $H\times H$ is regular subgroup of $G\wr S_2$. Suppose now that $G\in  \mathcal{K}_2(d)$.
Let $S\leq G$ be such that $S\in  
  \mathcal{K}_2(d)$.
    Since $S_\alpha\le G_\alpha$ and $S_\alpha\ne 1$, we have $|S_\alpha|=p$, implying $S_\alpha=G_\alpha$, and consequently $S=G$. This shows that no proper transitive subgroup of $G$ is in $K_2(d)$.
    By Lemma~\ref{lem:reduction}, there is $N\trianglelefteq G$ such that $|G/N|=2p^2$.
    Let $H=S_\alpha$ and let $A=HN/N$.
    
    Now we construct a complement of $A$ in $G/N$.
    If $p$ is odd, then $G/N$ has unique Sylow $p$-subgroup $P$ of order $p^2$.
    As $|A|=p$ and $f(A)\cap A=1$ for some $f\in\Aut(G/N)$, we have $P=Af(A)$, which implies that $P$ is not cyclic.
    Hence $P$ is elementary abelian. Let $y\in G/N$ be an involution.
    Then $G/N=P\rtimes\langle y\rangle\cong \mathbb{Z}_p^2\rtimes\mathbb{Z}_2$.
    As every irreducible representation of $\mathbb{Z}_2$ is 1-dimensional, there is a nontrivial $x\in P\setminus A$ such that $\mathbb{Z}_2$ fixes $\langle x\rangle$.
    Then $\langle x,y\rangle$ gives a complement of $A$, contradicting Lemma~\ref{lem:reduction}.
    
    If $p=2$, then $|G/N|=8$. As $A\cap f(A)=1$, $G/N$ has at least two involutions.
    This implies that $G/N\not\cong\mathbb{Z}_8$ and $G/N\not\cong Q_8$.
    In $\mathbb{Z}_2^3$, every subgroup has a complement.
    In $D_8$ and $\mathbb{Z}_2\times\mathbb{Z}_4$, the only subgroup of order $2$ without the complement is the characteristic subgroup generated by squares.
    As $A\cap f(A)=1$, $A$ is not characteristic, hence $A$ always has a complement, contradicting Lemma~\ref{lem:reduction}.
    This completes the proof.
\end{proof}

Combining Lemma~\ref{lem:reduction} with Corollary~\ref{cor:prime-stabilizer}, we obtain a useful
restriction on minimal members of $\mathcal K_2(d)$.
If $G\in\mathcal K_2(d)$ is minimal and $h=|G_\alpha|$,
then $h$ is a composite divisor of $d/2$.
Indeed, Lemma~\ref{lem:reduction} gives $|N|=d/(2h)$, so $h\mid d/2$, while
Corollary~\ref{cor:prime-stabilizer} excludes prime values of $h$.

\section{Groups of restricted degrees}\label{sec:small}
In this section, we continue our study of the wreath product $G \wr S_2$ when $G$ has restricted degree. 
By \Cref{thm:coprime-n-factorial}, if $G$ has odd degree, then $G \wr S_2$ in its product action contains a regular subgroup if and only if $G$ contains a regular subgroup. 
Therefore, we may assume that $G$ has even degree. It is then natural to consider groups of degree $2p,~ 4p$, or $8p$, where $p$ is a prime, as these form the simplest infinite families.

For the remainder of this section, $p$ is always assumed to be a prime.

\subsection{$G$ is of degree $2p$}
\begin{proposition}\label{prop:2p}
    Let $G\le\operatorname{Sym}(\Omega)$ be a transitive group of degree $2p$.
    Consider the product action of $G\wr S_2$ on $\Omega^2$.
    If $G\wr S_2$ has a regular subgroup on $\Omega^2$, then $G$ has a regular subgroup.
\end{proposition}
\begin{proof}
    It is equivalent to show that $\mathcal{K}_2(2p)$ is empty.
    Assume the contrary and choose a minimal counterexample $G$.
    Let $H$ be a point stabilizer of $G$.
   By Lemma~\ref{lem:reduction}, $|N|=p/|H|$, and since  $H\ne 1$, it follows that $|H|=p$, which contradicts to Corollary~\ref{cor:prime-stabilizer}.
\end{proof}

A straight consequence is that if a graph $\Gamma$ of order $2p$ is not isomorphic to a Cayley graph, then $\Gamma\mathbin{\square}\Gamma$ is not isomorphic to a Cayley graph.
This result also follows directly from Theorem~\ref{thm:main}.

\subsection{$G$ is of degree $4p$}

\begin{theorem}\label{thm:4p}
    Let $G\le\operatorname{Sym}(\Omega)$ be a transitive group of degree $4p$.
    If $G\wr S_2$ has a regular subgroup on $\Omega^2$, then $G$ has a regular subgroup.
\end{theorem}
\begin{proof}
    It is equivalent to show that $\mathcal{K}_2(4p)$ is empty.
    For $p\in\{2,3,7\}$, it can be verified computationally in \textsc{Magma} using $|G|\mid\frac{d^2}{2}$ and $|G_\alpha|\mid \frac{d}{2}$ by Lemma~\ref{lem:reduction}, so choose a minimal counterexample $G\in\mathcal{K}_2(4p)$, where $p\geq 5$ and $p\neq 7$.
    Then Lemma~\ref{lem:reduction} and Corollary~\ref{cor:prime-stabilizer} give $|G|=8p^2$, $|G_\alpha|=2p$, and $|N|=1$.
    Moreover, $B=\{(x,y)\in G\times G\mid y=f(x)\}$ for some $f\in\Aut(G)$.
    Let $P$ be a Sylow $p$-subgroup of $G$ and $Q$ a Sylow $2$-subgroup of $G$. Observe that 
 Corollary~\ref{cor:minimal-counterexample} implies that $P$ and $Q$ are not cyclic, as otherwise we obtain semiregular elements of order $p^2$ and $8$ in $G$ (respectively), which is impossible, as $G$ acts on $4p$ points.
 Hence $P\cong\mathbb{Z}_p^2$ and $Q$ has no element of order $8$.
 By a Sylow theorem, $P$ is normal in $G$.
    Write
    \[
 G=P\rtimes Q,\qquad G_\alpha=L\rtimes\langle b\rangle,
 \qquad |L|=p,\quad |b|=2,\quad b\in Q.
    \]
    If $Q$ has an elementary abelian subgroup $E$ such that $E\langle b\rangle=Q$, then $PE$ is transitive, supersolvable and every Sylow subgroup is elementary abelian.
    Then by \cite{Hall1937Comp}, $PE$ is complemented in $G$, and hence has a regular subgroup.
    This excludes $Q\cong\Z_{2}^3$ and also $Q\cong D_8$ when $b$ is a reflection.
    As $Q\not\cong\Z_2^{3}$, $Q$ has an element of order $4$.
    By Corollary~\ref{cor:minimal-counterexample}, after conjugation, $Q$ has a semiregular element $a$ of order $4$, which means $Q=\langle a\rangle\langle b\rangle$.
    If $Q\cong Q_8$ or $Q\cong D_8$ when $b$ lies in the center, then $Q=\langle a\rangle\langle b\rangle$ cannot be satisfied.
    The remaining case is $Q\cong\mathbb{Z}_2\times\mathbb{Z}_4$.
    In this case, $Q=\langle b\rangle\times\langle a\rangle$.

    If $\langle a\rangle$ acts reducibly on $P\cong\mathbb{Z}_p^2$, then $G$ has a normal subgroup $N$ of order $p$ such that $N\cap G_\alpha=1$. 
    Note that $\langle aN\rangle$ complements $G_\alpha N/N$ in $G/N$.
    Lemma~\ref{lem:regular-quotient-something} implies that $G$ has a regular subgroup, contrary to $G\in\mathcal{K}_2(4p)$.
    Hence $\langle a\rangle$ acts irreducibly on $P$.
    In particular, this action is faithful (since every involution in $\mathrm{GL}(2,p)$ is reducible), and therefore the minimal polynomial of $a$ on $P$ is $m(x)=x^2+1$, which implies that $a^2$ acts as $-I$.
    Since $b$ commutes with $a$, its $\pm1$-eigenspaces are $a$-invariant, so $b$ acts as $I$ or $-I$. For the former case, $\langle b\rangle$ is a nontrivial normal subgroup contained in $G_\alpha$, contradicting $\langle a\rangle$ being faithful on $P$. 
    Therefore $b$ acts as $-I$.

    The subgroup $P$ is characteristic in $G$. On $G/P\cong\Z_{2}\times\Z_{4}$, the automorphism induced by $f$ preserves the subgroup of squares $\langle Pa^2\rangle$. 
    Hence $f(Pb)\in\{Pb,Pa^2b\}$.
    Thus, $f(b)\in Pb\cup Pa^2b$.
    But $Pa^2b\subseteq C_G(P)$ and $b\notin C_G(P)$, so $f(b)\in Pb$.
    Now, for every $z\in P$,
    $b^z=z^{-1}bz=z^{-2}b.$        
    Since $p$ is odd, the map $z\mapsto z^{-2}$ is a bijection on $P$.
    Therefore every element of $Pb$ is conjugate to $b$ by an element of $P$.
    As $b\in G_\alpha$ fixes a point, every element of $Pb$ fixes a point as well.
    In particular, both $b$ and $f(b)$ fix points, and hence
    $(b,f(b))\in R$ fixes a point in the product action.
    This contradicts the fact that $B$ is semiregular.
\end{proof}

\subsection{$G$ is of degree $8p$}

In this subsection, we construct $G\in\mathcal{K}_2(8p)$ for every $p$ satisfying $p\equiv 3\pmod 4$.
Before we start, we need some elementary tools. For an element $g$ of a group $G$, we define $\Inn(g):G\to G$ by $(\Inn(g))(h)=h^g$ for every $h\in G$, the inner automorphism of $G$ corresponding to $g$.

\begin{lemma}
    \label{lem:when-contains-regular}
    Let $G$ be a group and let $H$ be a subgroup of $G$ such that $|G|=2|H|^2$.
    If $G$ has an automorphism $f$ and an element $g\in G$ such that 
    \begin{itemize}
        \item [\rm (i)] $f(H)\cap H=1$ and $f(H)\cap Hg=\emptyset$; and
        \item [\rm (ii)] $f^2=\operatorname{Inn}(g)$ and $f(g)=g$,
    \end{itemize}
    then $H\wr S_2$ has a complement in $G\wr S_2$.
\end{lemma}
\begin{proof}
    Let $S_2=\langle\sigma\rangle$ and let $h=|H|$.
    Then $G\wr S_2=(G\times G)\rtimes \langle \sigma\rangle$, where $\sigma$ interchanges the two copies of $G$.
    Let $B=\{(x,y)\in G\times G\mid f(x)=y\}$, let $\tau=(g,1)\sigma$
    and let $R=\langle B,\tau\rangle$.
    We claim that $R$ is a complement of $H\wr S_2$ in $G\wr S_2$.
    For $(x,y)\in B$, we have that $(x,y)^{\tau}=(y,x^g)=(f(x),f^2(x))\in B$ by (ii).
    Therefore, $\tau$ normalizes $B$.
    It follows that $R=B\langle \tau\rangle$.
    Since $\tau^2=(g,g)=(g,f(g))\in B$, we have that $|R|=2|B|$.
    Since $|B|=|G|=2h^2$, we have $|R|=2|B|=4h^2$.
    Now we claim that $G\wr S_2=R(H\wr S_2)$ and $R\cap (H\wr S_2)=1$.
    Let $s\in R\cap (H\wr S_2)$.
    Since $f(H)\cap H=1$ by (i), we have that $B\cap (H\times H)=1$.
    Therefore, if $s\in B$, then $s=1$.
    If $s\notin B$, then $s=(x,y)\tau$ for some $(x,y)\in B$.
    It follows that $s=(xg,y)\sigma\in (H\times H)\sigma$.
    Therefore, $xg\in H$ and $y\in H$.
    Since $f(x)=y$, we have that $yg=f(xg)\in Hg\cap f(H)=\emptyset$.
    This is a contradiction.
    Therefore, $s=1$ and $R\cap (H\wr S_2)=1$.
    Hence $$|R(H\wr S_2)|=4h^2\cdot |H\wr S_2|=|G\wr S_2|.$$
    This implies that $G\wr S_2=R(H\wr S_2)$.
\end{proof}

Now we illustrate that $\mathcal{K}_2(8p)$ is not empty when $p\equiv 3\pmod 4$.

\begin{example}\label{examp:8p}
    Let $p$ be a prime such that $p\equiv 3\pmod 4$ and
    let $G=\mathbb{Z}_2\times (\mathbb{Z}_p^2\times\mathbb{Z}_2^2)\rtimes\mathbb{Z}_4=\langle a\rangle\times (\langle x,y\rangle\times\langle b,c\rangle)\rtimes\langle d\rangle$,
    where $x^d=y^{-1}$, $y^d=x$, $b^d=c$ and $c^d=b$.
    Let $H=\langle x,ad^2,b\rangle$.
    We claim that $H$ is a core-free subgroup of $G$.
    Note that $H=\langle b\rangle\times\langle x\rangle\rtimes\langle ad^2\rangle\cong \mathbb{Z}_2\times\operatorname{D}_{2p}$.
    Then $H^d=\langle c\rangle\times\langle y\rangle\rtimes\langle ad^2\rangle$.
    It is clear that $H\cap H^d=\langle ad^2\rangle$ and $\langle ad^2\rangle$
    is not normal in $H$.
    Therefore, $\operatorname{Core}_G(H)< \langle ad^2\rangle$, which implies that $H$ is core-free in $G$.
    It follows that $G$ acts faithfully and transitively on the set $\Omega$ of right cosets of $H$ in $G$.
    
    We claim that $G$ does not have a regular subgroup.
Let $V=\langle x,y \rangle$ and let $T=\langle a,b,c,d^2,x,y \rangle$. Observe that $V$ is normal in $G$ and $T$ is of index $2$ in $G$. Suppose to the  contrary that $R$ is a regular subgroup of $G$. Then $|R|=8p$. Let $P=R\cap V$. Then $P$ has order $p$, and is normal in $R$. It follows that $R\leq N_G(P)$. Since $d^2$ acts as inversion on the elements of $V$, and $a,b,c,x,y$ act trivially on $V$, it follows that $T$ normalizes every subgroup of $V$, so in particular $T\leq N_G(P)$. Since $T$ is index 2 in $G$, $N_G(P)=T$ or $N_G(P)=G$. If $N_G(P)=G$, then $d$ normalizes $P$, hence it induces an automorphism of $P$ of order $4$, since $d^2$ acts as inversion. However, this is not possible, as $\Aut(P)$ has order $p-1\equiv 2 \pmod{4}$. We conclude that $N_G(P)=T$. Since both $R$ and $H$ are contained in $T$, it follows that $RH\subseteq T<G$, contradicting $G=RH$.

 Next, we show that $G\wr S_2$ has a regular subgroup.
    Let $S_2=\langle\sigma\rangle$.
    Then $G\wr S_2=(G\times G)\rtimes \langle \sigma\rangle$, where $\sigma$ interchanges the two copies of $G$.
    By quadratic reciprocity, there exists an integer $k$ such that $4k^4\equiv 1\pmod p$.
    Let $m$ be an integer such that $2mk\equiv 1\pmod p$.
    One can check that \[a\mapsto abc, b\mapsto ac, c\mapsto ab, d\mapsto d, x\mapsto x^ky^{-m}, y\mapsto x^my^{k}\]
    gives an automorphism $f$ of $G$.
    Moreover, one can check that $f$ and $g=d$ satisfy the conditions in Lemma~\ref{lem:when-contains-regular}.
    Therefore, $G\wr S_2$ has a regular subgroup.
\end{example}

    For $p=3$, the above construction recovers the permutation group of degree 24
(\texttt{TransitiveGroups(24)[675]}), which was the smallest-degree
example found in our \textsc{Magma} computations. 

\begin{remark}
If $G\wr S_2$ contains a regular subgroup $R$ in the product action,
then $G\wr S_{2k}$ contains a regular subgroup for every $k\geq1$.
By grouping the coordinates into $k$ pairs, we obtain a natural
embedding
\[
    (G\wr S_2)^k\leq G\wr S_{2k},
\]
under which $R^k$ acts regularly on $\Omega^{2k}$.
Thus Example~6.4 also provides examples for every even number
of factors. We are not aware of an analogous example with an
odd number of factors, which motivates the following question.
\end{remark}

\begin{problem}
Does there exist a finite transitive permutation group $G$
with no regular subgroup such that $G\wr S_n$ contains a
regular subgroup in the product action for some odd integer $n\geq3$?
\end{problem}

\begin{problem}
    Determine pairs $(d,n)$ of positive integers for which a transitive group $G$ of degree $d$ has a regular subgroup if and only if $G\wr S_n$ has a regular subgroup in the product action?
    Similarly, determine pairs $(d,n)$ for which every graph $\Gamma$ of order $d$ is Cayley if and only if $\Gamma^{\square n}$ is a Cayley graph?
\end{problem}

\section{A non-Cayley graph with Cayley Cartesian square}

Now we give a construction of a non-Cayley graph $\Gamma$ such that $\Gamma^{\mathbin{\square} 2}$ is isomorphic to a Cayley graph.
This construction is based on Example~\ref{examp:8p}.  Although the groups in that example have $2$-closures
containing regular subgroups, a modification of the
construction yields the required graphs.

Throughout this section, let $p \equiv 3 \pmod{4}$ be a prime,
let $\ell \geq 3$ be an odd integer coprime to $p$, and fix
$\rho \in \mathbb{Z}_\ell^*$ such that
$\rho^2 \equiv -1 \pmod{\ell}$.

\begin{example}\label{examp:8ml}
     Let $G$ be the group with generators $a,b,c,d,x,y$ defined in Example~\ref{examp:8p}, and let $G_{p,\ell}=\mathbb{Z}_\ell\rtimes G=\langle u\rangle\rtimes \langle a,b,c,d,x,y\rangle$, where $u^g=u$ for every $g\in\{a,b,c,x,y\}$ and $u^d=u^\rho$.
    Let $K_{p,\ell}=\langle b\rangle\times\langle x\rangle\rtimes\langle ad^2\rangle$.
    By Example~\ref{examp:8p}, $\langle u\rangle K_{p,\ell}/\langle u\rangle$ is core-free in $G_{p,\ell}/\langle u\rangle$.
    Then $\operatorname{Core}_{G_{p,\ell}}(K_{p,\ell})\le \operatorname{Core}_{G_{p,\ell}}(\langle u\rangle K_{p,\ell})\le \langle u\rangle$.
    Hence $\operatorname{Core}_{G_{p,\ell}}(K_{p,\ell})\le \langle u\rangle\cap K_{p,\ell}=1$.
    The action of $G_{p,\ell}$ on the coset space $[G_{p,\ell},K_{p,\ell}]$ is faithful.
    Hence, we treat $G_{p,\ell}$ as a transitive permutation group.
    Since $\gcd(p,\ell)=1$ and $\ell$ is odd,
    we have $\gcd(\ell,|G_{p,\ell}|/\ell)=1$.
    Also note that the action of $G_{p,\ell}/\langle u\rangle$ on the quotient space is permutation-isomorphic to the action of $G$.
    By Lemma~\ref{lem:regular-quotient-something}(iii), $G_{p,\ell}$ has no regular subgroup.
    Next, we show that $G_{p,\ell}\wr S_2$ in the product action has a regular subgroup. Let $W=G_{p,\ell}\wr S_2$.
    Note that $\langle u\rangle\times \langle u\rangle$ is a normal subgroup of $W$.
    Also note that the corresponding action of $W/(\langle u\rangle\times \langle u\rangle)$ and the action of $G\wr S_2$ are permutation isomorphic.
    Then by Lemma~\ref{lem:regular-quotient-something}(i), $W$ also has a regular subgroup. Therefore, $G_{p,\ell}\in\mathcal{K}_2(8p\ell)$.
\end{example}

\subsection{The graph construction}
\label{sec:graph}

We now construct a graph $\Gamma$ whose full automorphism group,
in its action on vertices, is permutation-isomorphic to $G_{p,\ell}$.
Denote $L=2p\ell$. Since $\gcd(2p,\ell)=1$, the Chinese remainder theorem gives a unique $\alpha\in\Z_L$ such that
\[ \alpha\equiv1\pmod 2, \qquad \alpha\equiv1\pmod p, \qquad \alpha\equiv\rho\pmod\ell. \]
The element $\alpha$ is a unit modulo $L$, since it is a unit modulo both $2p$ and $\ell$. Let
\[ q\colon\Z_L \to \Z_2\times\Z_\ell,\qquad x \mapsto (x\bmod2,x\bmod\ell), \]
and set
\[ S_0=\{(0,\pm 1),(1,\pm \rho)\},\qquad S=q^{-1}(S_0),\qquad K=\ker(q).\]
The homomorphism $q$ is surjective, so $|K|=p$. The four elements of $S_0$ are distinct, hence $|S|=4p$, $S+K=S$, and $S=-S$.

\begin{definition}\label{def:graph}
The graph $\Gamma=\Gamma_{p,\ell}^\rho$ has vertex set $V(\Gamma)=\Z_L \times \Z_2 \times \Z_2$, whose vertices are denoted by $(x,i,k)$. Its edge set is the union of the following three sets:
\begin{align*}
      R_1:\quad &(x,i,0)\sim(x\pm1,i,0),\\
      & (x,i,1)\sim(x\pm\alpha,i,1),\\
      R_2:\quad &(x,i,k)\sim(x,i+1,k),\\
      R_3:\quad &(x,i,0)\sim(z,j,1), \text{ for all } z\in x+S \text{ and } i,j\in\Z_2.
\end{align*}
\end{definition}

\begin{figure}[htbp]
\centering
\begin{tikzpicture}[
    scale=0.90,
    transform shape,
    every node/.style={font=\small},
    vertex/.style={circle,draw=black,fill=white,inner sep=1.7pt},
    Rone/.style={black,solid,line width=1pt},
    Rtwo/.style={black,dashed,line width=1pt},
    Rthree/.style={black,dotted,line width=1pt}
]
\def\rx{1.35}
\def\ry{0.50}
\def\rotang{18}
\def\upperY{0.68}
\def\lowerY{-0.68}
\def\layerX{3.65}

\coordinate (LTc) at (-\layerX,\upperY);
\coordinate (LBc) at (-\layerX,\lowerY);
\coordinate (RTc) at ( \layerX,\upperY);
\coordinate (RBc) at ( \layerX,\lowerY);

\foreach \i/\ang in {1/0,2/60,3/120,4/180,5/240,6/300}{
    \coordinate (LT\i) at
        ($(LTc)+({\rx*cos(\ang+\rotang)},{\ry*sin(\ang+\rotang)})$);
    \coordinate (LB\i) at
        ($(LBc)+({\rx*cos(\ang+\rotang)},{\ry*sin(\ang+\rotang)})$);
    \coordinate (RT\i) at
        ($(RTc)+({\rx*cos(\ang+\rotang)},{\ry*sin(\ang+\rotang)})$);
    \coordinate (RB\i) at
        ($(RBc)+({\rx*cos(\ang+\rotang)},{\ry*sin(\ang+\rotang)})$);
}

\foreach \i/\j in {1/2,2/3,3/4,4/5,5/6,6/1}{
    \draw[Rone] (LT\i) -- (LT\j);
    \draw[Rone] (LB\i) -- (LB\j);
    \draw[Rone] (RT\i) -- (RT\j);
    \draw[Rone] (RB\i) -- (RB\j);
}

\foreach \i in {1,...,6}{
    \draw[Rtwo] (LT\i) -- (LB\i);
    \draw[Rtwo] (RT\i) -- (RB\i);
}

\draw[Rthree] (LT2) to (RT3);
\draw[Rthree] (LT1) to (RB4);
\draw[Rthree] (LB1) to (RT4);
\draw[Rthree] (LB1) to (RB4);
\draw[Rthree] (LT1) to (RT4);

\foreach \i in {1,...,6}{
    \node[vertex] at (LT\i) {};
    \node[vertex] at (LB\i) {};
    \node[vertex] at (RT\i) {};
    \node[vertex] at (RB\i) {};
}

\node[font=\normalsize] at (-\layerX,1.45) {$k=0$};
\node[font=\normalsize] at ( \layerX,1.45) {$k=1$};

\node[left]  at (-5.35,\upperY) {$i=0$};
\node[left]  at (-5.35,\lowerY) {$i=1$};
\node[right] at ( 5.35,\upperY) {$i=0$};
\node[right] at ( 5.35,\lowerY) {$i=1$};

\begin{scope}[shift={(0,-1.80)}]
    \draw[Rone] (-5.80,0) -- (-5.00,0);
    \node[anchor=west] at (-4.85,0)
        {$R_1$: Cycle edges};

    \draw[Rtwo] (-1.90,0) -- (-1.10,0);
    \node[anchor=west] at (-0.95,0)
        {$R_2$: Matching edges};

    \draw[Rthree] (2.70,0) -- (3.50,0);
    \node[anchor=west] at (3.65,0)
        {$R_3$: Cross edges};
\end{scope}

\end{tikzpicture}
\caption{An illustration of the graph $\Gamma$. The number of vertices drawn on each cycle is for illustration only, and the $R_3$-edges shown are only representative.}
\label{fig:graph}
\end{figure}
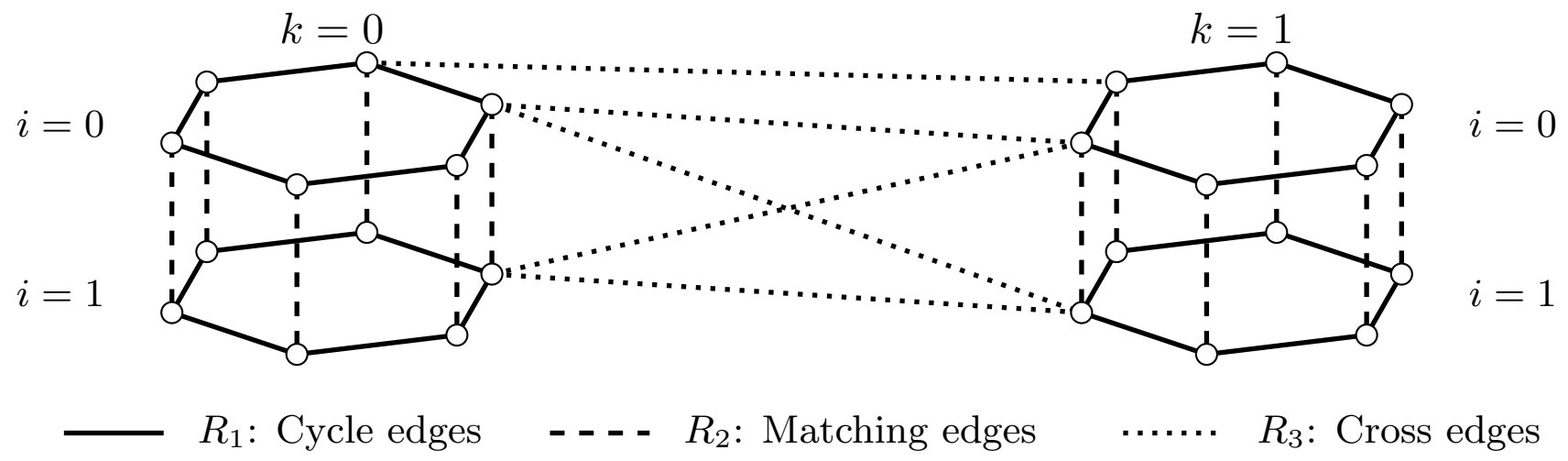

For $k\in\Z_2$, let
\[ V_k=\{(x,i,k):x\in\Z_L,\, i\in\Z_2\} \]
denote the \emph{$k$-layer} of $\Gamma$. Since both $1$ and $\alpha$ are units modulo $L$, the induced subgraph $\Gamma[V_k]$ is a prism $C_L\mathbin{\square} K_2$. The $R_3$-edges join the two prisms. Figure~\ref{fig:graph} illustrates this structure, where the cyclic orders in $V_0$ and $V_1$ are generated by addition by $1$ and by $\alpha$, respectively. 

For two vertices $(x,i,0)$ and $(z,j,1)$ in different layers, we define their \emph{cross difference} to be $z-x\in\Z_L$. They are adjacent if and only if their cross difference lies in $S$.

\begin{remark}
    This graph is obtained by experimental computation through \textsc{Magma}, and the smallest case is where $p=3$ and $\ell=5$ with 120 vertices.
    We computed the automorphism groups of generalized orbital graphs of $G_{3,5}$ and found this smallest example, which was then generalized to generalized orbital graphs of $G_{p,\ell}$.
\end{remark}

\begin{proposition}\label{prop:basic-graph}
      The graph $\Gamma$ is connected and simple, with $|V(\Gamma)|=8p\ell$ and valency $8p+3$.
\end{proposition}
\begin{proof}
    Each prism $\Gamma[V_k]$ is connected, and $R_3$ is nonempty, so $\Gamma$ is connected. Each vertex has two $R_1$-neighbors and one $R_2$-neighbor. It has $8p$ $R_3$-neighbors, corresponding to $|S|=4p$ choices of the first coordinate and $2$ choices of the second coordinate. Thus the valency is $8p+3$. The three edge sets are pairwise disjoint: An $R_1$-edge fixes both $i$ and $k$, an $R_2$-edge changes $i$ but fixes $k$, and an $R_3$-edge changes $k$. None of them contains a loop, so $\Gamma$ is simple. Finally, $|V(\Gamma)|=4L=8p\ell$.
\end{proof}

\subsection{The full automorphism group}

In this subsection, we determine the full automorphism group of $\Gamma$.

For $\varepsilon\in\{1,-1\}$, $a\in\Z_L$, $d\in K$, and $\eta_0,\eta_1\in\Z_2$, define
\begin{equation*}
P(\varepsilon,a,d,\eta_0,\eta_1)\colon
\begin{cases}
      (x,i,0)\longmapsto(\varepsilon x+a,i+\eta_0,0),\\
      (x,i,1)\longmapsto(\varepsilon x+a+d,i+\eta_1,1).
\end{cases}
\end{equation*}
Every such permutation preserves $R_1$ and $R_2$. Moreover, if an $R_3$-edge has cross difference $u\in S$, then its image has cross difference $\varepsilon u+d$, which also lies in $S$ because $S=-S$ and $S+K=S$. Thus every such permutation preserves $R_3$, and hence is an automorphism of $\Gamma$. The automorphisms $P(\varepsilon,a,d,\eta_0,\eta_1)$ are pairwise distinct and form a subgroup $A^+ \leq \Aut(\Gamma)$ of order $2\cdot L\cdot p\cdot4=16p^2\ell$.

Also define
\begin{equation*}
D\colon
\begin{cases}
      (x,i,0)\longmapsto(\alpha x,i,1),\\
      (x,i,1)\longmapsto(p\ell-\alpha^{-1}x,i,0).
\end{cases}
\end{equation*}
Since $\alpha$ is a unit modulo $L$, the map $D$ is a permutation of $V(\Gamma)$. It preserves $R_1$ and $R_2$; we now verify that it preserves $R_3$. Suppose that $(x,i,0)\sim(x+u,j,1)$ with $u\in S$, and write $q(u)=(\delta,\sigma\rho^\delta)$ with $\delta\in\Z_2$ and $\sigma\in\{1,-1\}$. After applying $D$ the new cross difference is
\begin{equation*}
      u'=\alpha x-p\ell+\alpha^{-1}(x+u) = (\alpha+\alpha^{-1})x+\alpha^{-1}u-p\ell.
\end{equation*}
It suffices to prove that $q(u')\in S_0$. Since $\rho^{-1}\equiv-\rho\pmod\ell$, we have $q(\alpha^{-1})=(1,-\rho)$ and $q(\alpha+\alpha^{-1})=(0,0)$. It follows that $q(u')=q(\alpha^{-1}u-p\ell)$. 
Moreover, $q(p\ell)=(1,0)$ and
\[ q(\alpha^{-1}u-p\ell) = (1,-\rho)(\delta,\sigma\rho^\delta)-(1,0) = (\delta+1,-\sigma\rho^{\delta+1})\in S_0. \]
It follows that $q(u')\in S_0$, hence $u'\in S$, so $D$ preserves $R_3$ and $D\in\Aut(\Gamma)$.

Define 
\[ A=\langle A^+,D\rangle\leq\Aut(\Gamma). \]
Observe that every element of $A^+$ preserves $V_0$ and $V_1$ setwise, whereas $D$ interchanges $V_0$ and $V_1$. 

\begin{lemma}\label{lem:intrinsic-edges}
      Every automorphism of $\Gamma$ preserves each of $R_1$, $R_2$ and $R_3$ setwise.
\end{lemma}
\begin{proof}
      For an edge $e$, let $\mu(e)$ denote the number of common neighbors of its endpoints. 

      Consider an $R_1$ edge in $V_0$. Using translations in $A^+$, we can assume that its endpoints are $(0,0,0)$ and $(1,0,0)$. The prism $\Gamma[V_0]$ is bipartite and contains no triangles. Therefore a common neighbor must lie in $V_1$, say $(z,t,1)$. Its first coordinate satisfies $z\in S$ and $z-1\in S$, while $t$ is arbitrary. Thus $q(z)$ belongs to $S_0\cap(S_0+(1,1))$. A direct computation shows that $|S_0\cap(S_0+(1,1))|$ equals $2$ if $\rho\equiv\pm2\pmod\ell$, and $0$ otherwise. Since $\rho^2\equiv-1\pmod\ell$, the first case occurs exactly when $\ell=5$. 
      Each point of $S_0$ has $p$ lifts to $S$, and $t$ has two choices. As $D$ interchanges the layers and preserves \(R_1\), we obtain, for every \(e\in R_1\), 
      \begin{equation*}
      \mu(e)=
      \begin{cases}
            4p,&\ell=5,\\
            0,&\ell\neq5.
      \end{cases}
      \end{equation*}

      The endpoints of an $R_2$-edge differ only in the second coordinate, while $R_3$-adjacency is independent of it. Thus they have exactly $8p$ common $R_3$-neighbors, and they have no common neighbor within their prism. Thus $\mu(e)=8p$ for every $R_2$ edge $e$. 
      Finally, let $e \in R_3$. The $R_2$-neighbor of each endpoint is adjacent to both endpoints, giving two common neighbors. Thus $\mu(e)\geq2$. On the other hand, any common neighbor in a fixed $k$-layer must be one of the two $R_1$-neighbors or the unique $R_2$-neighbor of the endpoint in that layer. Therefore $2\leq\mu(e)\leq6$, for every $R_3$ edge $e$.

      Since $p\geq 3$, the three edge sets have pairwise distinct values of $\mu$. As $\mu$ is invariant under graph automorphisms, every automorphism of $\Gamma$ preserves $R_1$, $R_2$ and $R_3$ setwise. 
\end{proof}

\begin{proposition}\label{prop:full-aut}
      $\Aut(\Gamma)=A$.
\end{proposition}
\begin{proof}
We already know that $A\leq \Aut(\Gamma)$, so it remains to prove the reverse inclusion.

      Deleting $R_3$ leaves the two connected prisms $\Gamma[V_0]$ and $\Gamma[V_1]$. By Lemma~\ref{lem:intrinsic-edges}, every graph automorphism of $\Gamma$ either preserves or interchanges them.
      Let $B$ be the subgroup of $\Aut(\Gamma)$ preserving both prisms setwise, then $|\Aut(\Gamma):B|\leq2$. Let $\varphi\in B$. As $\varphi$ preserves $R_1$ and $R_2$ separately, it has the form
      \begin{equation}\label{eq:Aplus-form}
      \begin{cases}
            (x,i,0)\longmapsto(\varepsilon_0x+a_0,i+\eta_0,0),\\
            (x,i,1)\longmapsto(\varepsilon_1x+a_1,i+\eta_1,1),
      \end{cases}
      \end{equation}
      where $\varepsilon_i\in\{1,-1\}$, $a_i\in\Z_L$, and $\eta_i\in\Z_2$. 
      Indeed, in either prism the two components of the $R_1$-subgraph are cycles.  The restriction to either cycle is dihedral, and the $R_2$-matching forces the same dihedral action on both cycles, apart from a possible simultaneous interchange of them.

      Set $d=a_1-a_0$. If an $R_3$-edge has cross difference $u\in S$, then its image under an automorphism of the form \eqref{eq:Aplus-form} has cross difference
      \begin{equation}\label{eq:image-cross-difference}
            (\varepsilon_1-\varepsilon_0)x+\varepsilon_1 u+d.
      \end{equation}
      If $\varepsilon_0\neq\varepsilon_1$, then $\varepsilon_1-\varepsilon_0=\pm2$. For fixed $u\in S$, the expression in \eqref{eq:image-cross-difference} runs through every residue modulo $\ell$ as $x$ varies, since $\ell$ is odd. This is impossible because the $\ell$-coordinate of an element of $S$ belongs to the proper subset $\{1,-1,\rho,-\rho\}$ of $\Z_\ell$. Therefore $\varepsilon_0=\varepsilon_1$ and we denote their common value by $\varepsilon$.
      Now \eqref{eq:image-cross-difference} reduces to $\varepsilon u+d$, which must belong to $S$ for all $u\in S$. Since $S=-S$, this is equivalent to $S+d=S$, or equivalently $S_0+q(d)=S_0$.

      We claim that $S_0+q(d)=S_0$ if and only if $q(d)=(0,0)$. Write $q(d)=(\delta,t)\in \Z_2\times\Z_\ell$. 
      Comparing the elements of $S_0$ and $q(d)+S_0$ with the first coordinate equal to $0$.
      When $\delta=0$, this implies that $\{\pm 1\}=\{\pm 1+t\}$, and when $\delta=1$, this implies $\{\pm 1\}=\{\pm\rho+t\}$.
      Since $\ell$ is odd, it follows that $t=0$. If $\delta=1$, comparing the elements of $S_0$ and $S_0+q(d)$ with first coordinate equal to 0 implies that $\{\pm 1\}=\{\pm\rho+t\}$. From here it follows that $2t=0$, but since $\ell$ is odd we conclude that $t=0$. However, we now obtain $\{\pm1\}=\{\pm\rho\}$, contrary to $\rho^2\equiv-1\pmod\ell$. Hence $q(d)=(0,0)$ and thus $d\in K$, and therefore
$\varphi=P(\varepsilon,a_0,d,\eta_0,\eta_1)\in A^+$. This shows that $B=A^+$.

Since $|\Aut(\Gamma):B|\leq 2$, and $D\not \in A^+$, it follows that $\Aut(\Gamma)=\langle A^+,D  \rangle=A.$
\end{proof}

\subsection{Cayley property of the graph and its Cartesian square}
In this subsection, we prove that $\Gamma$ is non-Cayley, whereas its Cartesian square $\Gamma \mathbin{\square} \Gamma$ is Cayley. To this end, we establish a permutation isomorphism between $A=\Aut(\Gamma)$ and the permutation group $G_{p,\ell}$ defined in Example~\ref{examp:8ml}.

Recall that $K=\ker(q)=\langle 2\ell\rangle$, and the following hold in $\Z_L$:
\[ \alpha p\ell =\alpha^{-1}p\ell=p \ell,\qquad \alpha\cdot 2\ell=\alpha^{-1}\cdot 2\ell=2\ell,\qquad \alpha\cdot 2p=-\alpha^{-1}\cdot 2p=\rho\cdot 2p. \]
Let $v_0=(0,0,0)\in V(\Gamma)$ be the base vertex. The vertex stabilizer of $v_0$ in $A$ is given by
\begin{equation}\label{eq:v0-stabilizer}
      A_{v_0} = \{ P(\varepsilon,0,\delta,0,\eta) \mid \varepsilon\in\{1,-1\},\, \delta\in K,\, \eta\in\Z_2 \}.
\end{equation}

\begin{lemma}
    \label{prop:perm-formula}
    The following statements hold.
    \begin{itemize}
        \item [\rm (i)] $P(1,a,d,\eta_0,\eta_1)\,P(1,a',d',\eta_0',\eta_1')=P(1,a+a',d+d',\eta_0+\eta_0',\eta_1+\eta_1')$.
        \item [\rm (ii)] $D^{-1}\,P(1,a,d,\eta_0,\eta_1)\,D=P(1,-\alpha^{-1}(a+d),\alpha a+\alpha^{-1}(a+d),\eta_1,\eta_0)$.
    \end{itemize}
\end{lemma}
\begin{proof}
     {\rm (i)} This follows directly from the definition of $P$.

    {\rm (ii)} Let $P=P(1,a,d,\eta_0,\eta_1)$. Recall that we use the right-action convention. Thus the product $D^{-1}PD$ is evaluated by applying $D^{-1}$, then $P$, and finally $D$. We first compute $D^{-1}$. Since $\alpha\equiv1\pmod{2p}$ and $L=2p\ell$, we have $\alpha p\ell=\alpha^{-1}p\ell=p\ell$ in $\Z_L$. Therefore
    \[ D^2(x,i,k)=(p\ell-x,i,k), \]
    so $D^4=1$ and $D^{-1}=D^3$. Applying $D$ once more gives
    \[ D^{-1}(x,i,0)=(p\ell-\alpha x,i,1),\qquad D^{-1}(x,i,1)=(\alpha^{-1}x,i,0). \]
    Now let $v=(x,i,k)\in V(\Gamma)$. If $k=0$, then
    \begin{align*}
        v^{D^{-1}}&=(p\ell-\alpha x,i,1),\\
        \bigl(v^{D^{-1}}\bigr)^{P}&=(p\ell-\alpha x+a+d,\,i+\eta_1,\,1),\\
        \bigl(v^{D^{-1}P}\bigr)^{D}&=\bigl(p\ell-\alpha^{-1}(p\ell-\alpha x+a+d),\,i+\eta_1,\,0\bigr)
        =\bigl(x-\alpha^{-1}(a+d),\,i+\eta_1,\,0\bigr),
    \end{align*}
    where the last equality uses $\alpha^{-1}p\ell=p\ell$. If $k=1$, then
    \[ v^{D^{-1}}=(\alpha^{-1}x,i,0),\\
       \bigl(v^{D^{-1}}\bigr)^{P}=(\alpha^{-1}x+a,i+\eta_0,0),\\
       \bigl(v^{D^{-1}P}\bigr)^{D}=(x+\alpha a,i+\eta_0,1). \]
    Let
    \[ a'=-\alpha^{-1}(a+d),\qquad d'=\alpha a+\alpha^{-1}(a+d). \]
    It remains to check that $d'\in K$. Since $q(\alpha)=(1,\rho)$ and $q(\alpha^{-1})=(1,\rho^{-1})=(1,-\rho)$, we have $q(\alpha+\alpha^{-1})=(0,0)$, that is, $\alpha+\alpha^{-1}\in K$. 
    Since $\alpha+\alpha^{-1}\equiv0\pmod{2\ell}$ and
$d\equiv0\pmod{2\ell}$, we have
\[
    d'=(\alpha+\alpha^{-1})a+\alpha^{-1}d
       \equiv0\pmod{2\ell}.
\]
Thus $d'\in K$, completing the proof.
\end{proof}

\begin{proposition}\label{prop:permutation-isomorphism}
      The map $\Phi:G_{p,\ell}\longrightarrow A$ defined on the generators of $G_{p,\ell}$ by
      \begin{equation*}
      \begin{aligned}
            a&\longmapsto \ha=P(1,p\ell,0,0,0),&
            b&\longmapsto \hb=P(1,0,0,0,1),\\
            c&\longmapsto \hc=P(1,0,0,1,0),&
            x&\longmapsto \hx=P(1,0,2\ell,0,0),\\
            y&\longmapsto \hy=P(1,2\ell,-2\ell,0,0),&
            u&\longmapsto \hu=P(1,2p,0,0,0),\\
            d&\longmapsto \hd=D.
      \end{aligned}
      \end{equation*}
      is a group isomorphism. Moreover, the actions of $G_{p,\ell}$ in right coset action and $A$ in $V(\Gamma)$ are permutation isomorphic.
\end{proposition}
\begin{proof}
      By Lemma~\ref{prop:perm-formula}, $|\ha|=|\hb|=|\hc|=2$, $|\hd|=4$, $|\hu|=\ell$, and $|\hx|=|\hy|=p$.
      First, we claim that $A$ is generated by $\ha,\hb,\hc,\hx,\hy,\hu,\hd$.
      By Lemma~\ref{prop:perm-formula}~{\rm(i)}, $\langle \ha,\hb,\hc,\hx,\hy,\hu\rangle$ is abelian.
      Denote this subgroup by $T$.
      One can check that $|T|=2^3\cdot p^2\cdot\ell=8p^2\ell$ and $\langle\hd\rangle\cap T=1$.
      Therefore, $|\langle \ha,\hb,\hc,\hx,\hy,\hu,\hd\rangle|\ge |T||\hd|=32p^2\ell=|A|$. This implies that $A$ is generated by $\ha,\hb,\hc,\hx,\hy,\hu,\hd$.
      By Lemma~\ref{prop:perm-formula} {\rm (ii)}, one can check that
      $\ha^{\hd}=\ha$, $\ \hb^{\hd}=\hc$, $\ \hc^{\hd}=\hb$, $\ \hx^{\hd}=\hy^{-1}$, $\ \hy^{\hd}=\hx$, $\ \hu^{\hd}=\hu^{\rho}$.
      Therefore, the map $\Phi$ is a well-defined homomorphism and is surjective, as $A$ is generated by $\ha,\hb,\hc,\hx,\hy,\hu,\hd$.
      Since $|G_{p,\ell}|=|A|$, we have that $\Phi$ is bijective, and hence is a group isomorphism.

      To complete the proof, it remains to show that $\Phi$ induces a permutation isomorphism between the two actions. The group $A^+$ is transitive on each layer, and $D$
interchanges the layers, so $A$ is transitive on $V(\Gamma)$. Since both actions are transitive, by \cite[Lemma 2.6]{Praeger_Schneider_2018}, it suffices to show that $K_{p,\ell}^\Phi=A_{v_0}$, where $K_{p,\ell}=\langle x,b,ad^2 \rangle$ is a stabilizer of $G_{p,\ell}$ in the right coset action; see Example~\ref{examp:8ml}. By the definition of $\Phi$,
      \[ K_{p,\ell}^\Phi = \langle x,b,ad^2 \rangle^\Phi=\langle \hx,\hb,\ha D^2 \rangle. \]
      On the other hand, the stabilizer in \eqref{eq:v0-stabilizer} gives $A_{v_0}=\langle \hx,\hb,\ha D^2 \rangle$. Therefore $K_{p,\ell}^\Phi=A_{v_0}$, and the required permutation isomorphism follows.
\end{proof}

\begin{proof}[Proof of Theorem~\ref{thm:intro-counterexamples}]
      By Proposition~\ref{prop:full-aut}, we have $A=\Aut(\Gamma)$. Moreover, Proposition~\ref{prop:permutation-isomorphism} shows that $A$ is permutation isomorphic to $G_{p,\ell}$. By Example~\ref{examp:8ml}, the group $G_{p,\ell}$ contains no regular subgroup, while $G_{p,\ell}\wr S_2$ contains a regular subgroup. Therefore, $\Gamma$ is non-Cayley, whereas $\Gamma \mathbin{\square} \Gamma$ is Cayley.
\end{proof}

The group and graph versions of the problem can behave differently
at a fixed degree. Indeed, Example~\ref{examp:8p} shows that
$\mathcal K_2(24)\neq\emptyset$. However, using the complete list
of connected vertex-transitive non-Cayley graphs of order $24$ \cite{HOLT2020},
we verified in \textsc{Magma} that none of their Cartesian squares is
Cayley. 
Example~\ref{examp:8p} establishes that
$\mathcal K_2(8p)\neq\emptyset$ for every prime
$p\equiv3\pmod4$. This raises the following question.

\begin{problem}
For which primes $p\equiv3\pmod4$ does there exist a connected
non-Cayley graph of order $8p$ whose Cartesian square is Cayley?
\end{problem}

\begin{remark}
    The graph $\Gamma$ also provides an example of a non-Cayley graph such that $\Gamma \times \Gamma$ (tensor or direct product) and $\Gamma \boxtimes \Gamma$ (strong product) are Cayley graphs, as the regular subgroup of the wreath product constructed in Example~\ref{examp:8ml} also acts regularly on both products.
\end{remark}

\section*{Acknowledgment}
Zhaochen Ding is supported by the Slovenian Research and Innovation Agency (ARIS) through research program P1-0285.
ARIS supports Ted Dobson through research program P1-0285 and research projects J1-70047, J1-70035, J1-70046 and J1-50000.
ARIS supports Ademir Hujdurović through research program P1-0404 and research projects N1-0481, N1-0428, N1-0459, J1-70035, J1-70046, J1-70047, and J1-60012. 

\section*{Declaration of AI use}
AI was used for minor polishing of the presentation; the authors take full responsibility for all the results.

\bibliography{References}{}
\bibliographystyle{amsplain}
\end{document}